\documentclass[12pt, a4paper]{article}
\usepackage[utf8]{inputenc}
\usepackage{amssymb,amsthm,amsmath}
\title{Cohomological properties of Hermitian sympletic threefolds}
\author{Grigory Papayanov}
\date{2015}

\def\blacksquare{\hbox{\vrule width 7pt height 7pt depth 0pt}}
\def\endproof{\blacksquare}

\newcommand{\one}{\omega^{1,1}}

\newcommand{\arrow}{{\:\longrightarrow\:}}
\newcommand{\Z}{{\Bbb Z}}
\newcommand{\C}{{\Bbb C}}
\newcommand{\R}{{\Bbb R}}

\renewcommand{\H}{{\mathcal{H}}}
\newcommand{\I}{\mathcal{I}}
\newcommand{\6}{\partial}
\def\1{\sqrt{-1}\:}
\renewcommand{\bar}{\overline}
\renewcommand{\phi}{\varphi}
\renewcommand{\epsilon}{\varepsilon}

\newcommand{\Vol}{\operatorname{Vol}}

\newcommand{\A}{\mathcal{A}}
\newcommand{\Ker}{\operatorname{Ker}}
\newcommand{\ddc}{dd^c}
\newcommand{\Alb}{\operatorname{Alb}}
\renewcommand{\dim}{\operatorname{dim}}

\renewcommand{\Im}{\operatorname{Im}}

\newcounter{Mycounter}[section]
\newcounter{lemma}[section]
\renewcommand{\thelemma}{\noindent{Lemma \thesection.\arabic{lemma}}}
\newcommand{\lemma}{%
     \setcounter{lemma}{\value{Mycounter}}
     \refstepcounter{lemma}
     \stepcounter{Mycounter}
     {\bf \thelemma:\ }}

\newcounter{claim}[section]
\newcounter{sublemma}[section]
\newcounter{corollary}[section]
\renewcommand{\thecorollary}{\noindent{Corollary \thesection.\arabic{corollary}}}
\newcommand{\corollary}{%
     \setcounter{corollary}{\value{Mycounter}}
     \refstepcounter{corollary}
     \stepcounter{Mycounter}
     {\bf \thecorollary:\ }}

\newcounter{theorem}[section]
\renewcommand{\thetheorem}{\noindent{Theorem \thesection.\arabic{theorem}}}
\newcommand{\theorem}{%
     \setcounter{theorem}{\value{Mycounter}}
     \refstepcounter{theorem}
     \stepcounter{Mycounter}
     {\bf \thetheorem:\ }}

\newcounter{conjecture}[section]
\newcounter{proposition}[section]
\renewcommand{\theproposition}
       {\noindent{Proposition \thesection.\arabic{proposition}}}
\newcommand{\proposition}{%
     \setcounter{proposition}{\value{Mycounter}}
     \refstepcounter{proposition}
     \stepcounter{Mycounter}
     {\bf \theproposition:\ }}

\newcounter{definition}[section]
\renewcommand{\thedefinition}
       {\noindent{Definition~\thesection.\arabic{definition}}}
\newcommand{\definition}{%
     \setcounter{definition}{\value{Mycounter}}
     \refstepcounter{definition}
     \stepcounter{Mycounter}
     {\bf \thedefinition:\ }}

\newcounter{example}[section]
\newcounter{remark}[section]
\renewcommand{\theremark}{\noindent{Remark \thesection.\arabic{remark}}}
\newcommand{\remark}{%
     \setcounter{remark}{\value{Mycounter}}
     \refstepcounter{remark}
     \stepcounter{Mycounter}
     {\bf \theremark:\ }}

\newcounter{problem}[section]
\newcounter{question}[section]
\begin{document}

\begin{center}
{\Large\bf
Cohomological properties of Hermitian\\ symplectic threefolds
}
%%%%%%%%%%%%%%%%%%%%%%%%%%%%%%%%%%%%%%%%%%%%%%%%%%%%%%%%%%%%
\\[4mm]

Grigory Papayanov
\\[6mm]

\end{center}

\begin{abstract}

A Hermitian symplectic manifold is a complex manifold endowed with a symplectic form
$\omega$, for which the bilinear form $\omega(I\cdot,\cdot)$ is positive definite. In this work
we prove $\ddc$-lemma for 1- and (1,1)-forms for compact Hermitian symplectic 
manifolds of dimension~3. This shows that Albanese map for such manifolds is well-defined
and allows one to prove K\"ahlerness if the dimension of the Albanese image of a manifold is maximal.
\end{abstract}

%\renewcommand{\cftsecleader}{\cftdotfill{\cftdotsep}}
%\tableofcontents

\section*{Introduction}
\addcontentsline{toc}{section}{Introduction}

A Hermitian symplectic manifold is a complex manifold $(M,I)$ together with a symplectic form $\omega$,
 for which the bilinear form
$\omega(I\cdot,\cdot)$ is positive definite (that is, $\omega(IX,X)>0$ for any vector field $X$ on $M$).
Any K\"ahler manifold is obviously Hermitian symplectic,
and it is an open problem whether there exist other examples
of Hermitian symplectic manifolds. Hermitian symplectic
manifolds were studied by Streets and Tian in 
\cite{Streets_Tian:pluriclosed} and 
\cite{Streets_Tian:flow}; they constructed an appropriate Ricci
flow on Hermitian symplectic
manifolds, and studied its convergency
properties. Since then, many people searched for non-trivial 
examples of Hermitian symplectic manifolds.

The search for non-K\"ahler examples of Hermitian symplectic
manifolds was vigorous, but ultimately unsuccessful.
All common sources of examples of non-K\"ahler manifolds
were tapped at some point. 

For complex dimension 2, Hermitian symplectic structures
are all K\"ahler. This was shown by Streets and Tian in
\cite{Streets_Tian:pluriclosed}. Another proof could be
obtained from the Lamari (\cite{Lamari}) result about
existence of positive, exact $(1,1)$-current on any
non-K\"ahler complex surface.

In \cite{Peternell}, it was shown that any 
non-K\"ahler Moishezon manifold admits an exact, positive
$(n-1,n-1)$-current; therefore, Moishezon manifolds
which are Hermitian symplectic are also K\"ahler.

In \cite{Enrietti_Fino_Vezzoni} it was shown 
that no complex nilmanifold can admit a Hermitian
symplectic structure, and in 
\cite{Fino_Kasuya_Vezzoni} this result was extended 
to all complex solvmanifolds and Oeljeklaus-Toma manifolds.

Existence of K\"ahler metric implies some restrictions on the
cohomology of a manifold: for example the Fr\"olicher spectral
sequence of K\"ahler manifold always degenerates at the first page.
Results of Cavalcanti (\cite{Cavalcanti:SKT}) show that
the Fr\"olicher spectral sequence for Hermitian symplectic manifolds degenerates at the first page.

In this work we define some Laplacian-like operators, kernels of which
conjecturally isomorphic to the spaces of cohomology, and, with the help of these operators,
prove $\ddc$-lemma for (1,1)-forms on Hermitian symplectic threefolds. Argument of Gauduchon
(\cite{Gauduchon}) shows that $\ddc$-lemma for (1,1)-forms is equivalent to the equality
$b^1=2h^{0,1}$. It follows that the Albanese map is well-defined and, if its image is not a point,
the generic fiber of $\Alb$ is K\"ahler. The question of existence of special (e.g. K\"ahler or balanced)
metrics on total spaces of maps with K\"ahler base and fibers is studied, for example, in \cite{HL}
and \cite{Michelsohn}. Using the Albanese map, we are able to prove that if a Hermitian symplectic threefold $M$
has $\dim \Alb(M)=3$, then it admits a K\"ahler metric, 
and if $\dim \Alb(M)=1$, $M$ is balanced. If $dd^c$-lemma holds for $(2,2)$-forms, 
then by \cite{HL} $\dim \Alb(M)=2$ would imply that $M$ is K\"ahler, 
but, unfortunately, we have not proven $\ddc$-lemma in full generality yet.

{\bf Acknowledgements.} The author would like to thank M.\,Verbitsky for many
extremely helpful discussions. Work on sections 1--3 was supported by RSCF,
 grant number 14-21-00053, within the
Laboratory of Algebraic Geometry. Work on section 4 was supported by RFBR 15-01-09242.

%\section*{Keywords}

%Hermitian symplectic manifolds, SKT manifolds, Fr\"olicher spectral sequence, $dd^c$-lemma.

%%%%%%%%%%%%%%%%%%%%%%%%%%%%%%%%%%%%%%%%%%%%%%%%%%%%%%%%%%%%%%%%%%%%%%%
\section{Preliminaries}
%%%%%%%%%%%%%%%%%%%%%%%%%%%%%%%%%%%%%%%%%%%%%%%%%%%%%%%%%%%%%%%%%%%%%%%

\definition
Let $M$ be a smooth manifold of dimension 2n, $I:TM \arrow TM$ an integrable
complex structure, $\A^{p,q}$ the corresponding Hodge decomposition on the
bundle of differential forms: $\A^n\otimes \C=\bigoplus_{n=p+q}\A^{p,q}$, $\one$ a form
in $\A^{1,1}$. We will say that $\one$ is {\it Hermitian} if the tensor
$h(\cdot,\cdot):=\one(\cdot,I\cdot)$ is a 
Riemannian metric on $M$, and we will say that $\one$ is {\it Hermitian symplectic}
if there exists a symplectic form $\omega$ such that $\one$ is the (1,1)-component
in the Hodge decomposition of $\omega$. If $M$ is endowed with such $\I$ and $\one$, we will
call it a Hermitian symplectic manifold.

\hfill

For a Hermitian symplectic manifold $(M,I,\omega)$, let $d: \A^\bullet\arrow \A^{\bullet+1}$ be the usual de Rham
differential acting on forms, $d^c:=IdI^{-1}: \A^{\bullet}\arrow \A^{\bullet+1}$  the twisted differential,
$L: A^\bullet\arrow A^{\bullet+2}$  the operator of (left) multiplication by $\omega$,
$L(\eta):= \omega\wedge \eta$, $\Lambda: \A^{\bullet}\arrow \A^{\bullet-2}$ the adjoint operator (\cite{Yau_Tseng}).
In the local Darboux coordinates $p_i, q_i$ where $\omega=\sum dp_i\wedge dq_i$, operator $\Lambda$  
looks like $\sum i_{\!\frac{\6}{\6p_i}}i_{\!\frac{\6}{\6q_i}}$. We will denote by $L^{1,1}$ the
operator of multiplication by the hermitian form $\one$, and by $\Lambda^{1,1}$ the adjoint operator to $L^{1,1}$.

\hfill

\lemma\label{SKT}
The form $\one$ is the SKT form, that is, $\6\bar\6\one=0$.

\hfill

{\bf Proof:}
Let $\omega=\one+\alpha$, where $\alpha$ lies in $\A^{2,0}\oplus\A^{0,2}$.
Since $d\omega=0$, $\6\one=-\bar\6\alpha$ and $\6\bar\6\one=\bar\6^2\alpha=0$. \endproof

%\remark
%For the information about Hermitian symplectic manifolds see
%Streets-Tian, Popovici, Fino etc. (REFERENCES).

\hfill

\definition
Let $\alpha$ be a differential form on $M$. We will say that $\alpha$
is {\it primitive with respect to $\omega$} if $\Lambda\alpha=0$, and that
$\alpha$ is primitive with respect to $\one$ if $\Lambda^{1,1}\alpha=0$.

\hfill

\lemma (The Weil identities).
Let $B^{p,q}$ be a primitive with respect to $\one$ $(p,q)$-form, $p+q=r$. Then
the following formula holds (\cite[Proposition 6.29]{Voisin}):

$$*B^{p,q}=(-1)^{\frac{r(r+1)}{2}}(\sqrt{-1})^{p-q}\frac{1}{(n-r)!}(\one)^{n-k}\wedge B^{p,q}.$$

\endproof

\hfill

\definition
An operator $\Delta$ defined as double graded commutator,
$\Delta:=\{d,\{d^c,\Lambda^{1,1}\}\}$
is called {\it the Hermitian symplectic} Laplacian.

\hfill

\remark
$\Delta$ is not a Laplacian associated to the Riemannian metric $h$.
Nevertheless they differ by a differential operator of
first order (see e.g. \cite{Liu_Yang} for the exact formula),
 therefore they have equal symbols, so $\Delta$ is elliptic.

\hfill

Recall the graded Jacobi identity for the graded commutator:
$$\{a,\{b,c\}\}=\{\{a,b\},c\}+(-1)^{deg(a)deg(b)}\{b,\{a,c\}\}. $$

\hfill

\lemma\label{commutators}
$\Delta=\{d^c,\{d,\Lambda^{1,1}\}\}$. Therefore $\Delta$ commutes 
with $d$ and with $d^c$.

\hfill

{\bf Proof:} Follows simply from the Jacobi identity.
\endproof

\hfill

\theorem\label{spectral} (Spectral theorem).
Let $(M, I,\omega)$ be a compact Hermitian symplectic manifold. Then
the space of differential forms decomposes as a topological direct sum of
generalized eigenspaces of $\Delta$:
 $\A^\bullet(M)=\bigoplus_{\lambda_i}\A^\bullet_{\lambda_i}(M)$,
each component of this decomposition is finite-dimensional and
 preserved by $d$, $d^c$ and $\delta$.

\hfill

{\bf Proof:}
Decomposition is in fact proven in \cite[Proposition 2.36]{BGV}
($\Delta$ is a generalized laplacian in their terminology);
one has
to apply spectral theorem for compact operators: compact operator on
Hilbert space has a canonical Jordan form with finite-dimensional generalized
eigenvalues (\cite{Conway}).

By \ref{commutators}, $\Delta$ commutes with $d$ and $d^c$, so all
generalized eigenspaces are in fact subcomplexes.
\endproof

\hfill

\theorem
Let $\alpha$ be a closed form in $\bigoplus_{\lambda_i\ne 0}\A^\bullet_{\lambda_i}(M)$. Then
$\alpha$ is exact.

\hfill

{\bf Proof:}
When restricted to $\bigoplus_{\lambda_i\ne 0}\A^\bullet_{\lambda_i}(M)$, Laplacian
$\Delta$ has an inverse, $\Delta^{-1}$. So
$$\alpha=\Delta\Delta^{-1}\alpha=(\pm dd^c\Lambda \pm d\Lambda d^c)\Delta^{-1}\alpha.$$
\endproof

\hfill

%%%%%%%%%%%%%%%%%%%%%%%%%%%%%%%%%%%%%%%%%%%%%%%%%%%%%%%%%%%%%%%%%%%%%%%
\section{Forms on a Hermitian symplectic manifold}
%%%%%%%%%%%%%%%%%%%%%%%%%%%%%%%%%%%%%%%%%%%%%%%%%%%%%%%%%%%%%%%%%%%%%%%

In this section $M$ is assumed to be compact.

\hfill

\lemma\label{ddc1} ($dd^c$-lemma for 1-forms).
Let $\alpha$ be a $d$-exact, $d^c$-closed (or $d^c$-exact and $d$-closed)
1-form. Then $\alpha=0$.

\hfill

{\bf Proof:}
Suppose $\alpha$ is $d$-exact, $\alpha=df$. Then $dd^cf=0$. By Hopf maximum principle
(\cite{_Gilbarg_Trudinger_}), $f$ is constant, hence $\alpha=df=0$. \endproof

\hfill

We will now investigate whether holomorphic forms on $M$ are closed.

\lemma
Let the $n$ be the complex dimension of $M$. Then every holomorphic $n-2$-form is closed.

\hfill

{\bf Proof:} Let $\alpha$ be a holomorphic $n-2$-form, $\alpha \in \A^{n-2,0}$,
$\bar\6\alpha=0$. Then $d\alpha=\6\alpha$ is primitive with respec to $\one$,
by dimension reasons. So, by Weil identities,
$$||d\alpha||^2=\int d\alpha\wedge d\bar\alpha \wedge \one
=\int \6\alpha\wedge \bar\6\bar\alpha \wedge \one = 
\alpha \wedge \bar\alpha \wedge \6\bar\6\one=0.$$
Hence $\alpha$ is closed. \endproof

\hfill

\remark\label{holoforms}
Obviously, on any compact complex manifold of complex dimension $n$, 
every holomorphic function and every holomorphic $n$-form is closed.
Every holomorphic $n-1$-form is also closed, as the simple argument with the
integration shows. So, any holomorphic form on a Hermitian symplectic 
threefold is closed.

\hfill

%%%%%%%%%%%%%%%%%%%%%%%%%%%%%%%%%%%%%%%%%%%%%%%%%%%%%%%%%%%%%%%%%%%%%%%
\section{$dd^c$-lemma for (1,1)-forms}
%%%%%%%%%%%%%%%%%%%%%%%%%%%%%%%%%%%%%%%%%%%%%%%%%%%%%%%%%%%%%%%%%%%%%%%

Recall that by \ref{spectral} every differential form $\alpha$ decomposes
by generalized eigenspaces of $\Delta$: $\alpha=\alpha_0 + \alpha_{\ne 0}$,
where $\Delta^N(\alpha_0)=0$ for some $N$, and 
$\alpha_{\ne 0}=\Delta\Delta^{-1}\alpha_{\ne 0}$. Suppose that $\alpha$ is
$d$-exact and $d^c$-closed. Then $\alpha_0$ and $\alpha_{\ne 0}$ are also
$d$-exact and $d^c$-closed.

\hfill

\lemma
In notations as above, $\alpha_{\ne 0}$ is $dd^c$-exact.

\hfill

{\bf Proof:} by \ref{commutators}, $\Delta^{-1}$ commutes with $d$ and $d^c$,
so $\Delta\Delta^{-1}\alpha_{\ne 0}=dd^c\Lambda^{1,1}\Delta^{-1}\alpha_{\ne 0}=\alpha$.
\endproof

\hfill

\lemma\label{primitivness}
Suppose exact (1,1)-form $\eta=d\gamma$ lies in the kernel of $\Delta^{1,1}$. Then 
$\eta$ is primitive (with respect both to $\omega$ and to $\omega^{1,1}$). 

\hfill

{\bf Proof:}
$\Delta\eta=dd^c\Lambda^{1,1}\eta=0$, so, by Hopf maximum principle \cite{_Gilbarg_Trudinger_}
$\Lambda^{1,1}\eta=c$, where $c$ is some constant. It means that $\Lambda\eta$ also equals $c$.

If $\Lambda\eta=c,$ then $\eta=c\omega+B$, where $B$ is a primitive form. Since $\eta=d\gamma$,
the cohomology classes of $c\omega$ and $B$ are equal, but the cohomology class of a symplectic
form cannot be represented by a primitive form. Indeed, $\omega\wedge\omega^{n-1}$ is a volume form,
hence nonzero in cohomology, but $B\wedge\omega^{n-1}=0$. So $c=0$ and $\eta$ is primitive.
\endproof

\hfill

\lemma\label{vanish}
Suppose $\dim(M)=3$, $\eta=dd^cf$ is $dd^c$-exact primitive $(1,1)$-form. Then $\eta=0$.

\hfill

{\bf Proof:}
Note first that, since $\eta$ is primitive with respect to $\omega$, it is primitive
with respect to $\omega^{1,1}$, so, by Weil identities,
$*\eta=\eta\wedge(\omega^{1,1})^{\wedge n-2}$, where
$*$ is the Hodge star operator associated with the Hermitian metric $h$ with corresponding
2-form equal to $\omega^{1,1}$ (\cite{Griffiths_Harris}). Then $h(\eta,\eta)=$

$$
=\int \eta\wedge*\eta=\int \eta\wedge\eta\wedge(\omega^{1,1})^{\wedge n-2}=\int f\eta\wedge dd^c(\omega^{1,1})^{\wedge n-2}.
$$

But $dd^c\omega^{1,1}=0$ on a Hermitian symplectic manifold, so the integral vanishes.
Since $h$ is a hermitian metric, $\eta$ also equals to zero.
\endproof

\hfill

\lemma
Let $\dim(M)=3$. Suppose that an exact (1,1)-form $\eta=d\gamma$ lies in
the kernel of $(\Delta)^n,$ $n>1$. Then 
$\eta$ lies in the kernel of $(\Delta)^{n-1}.$

\hfill

{\bf Proof:}
$(\Delta)^{n-1}\eta$ is an exact (1,1)-form lying in the kernel of $\Delta$, so,
by \ref{primitivness} it is primitive. Since $d\eta=d^c\eta=0, (\Delta)^{n-1}\eta=(dd^c\Lambda)^{n-1}\eta$,
it is $dd^c$-exact, therefore, by \ref{vanish}, it vanishes.
\endproof

\hfill

In order to complete the proof of $\ddc$-lemma for (1,1)-forms on Hermitian symplectic manifolds, we have to 
prove that an exact, primitive (1,1)-form vanishes.

\hfill

\lemma 
Let $M$ be a Hermitian symplectic manifold of dimension 3, $\eta$ be an exact, primitive
(1,1)-form on $M$. Then $\eta=0$.

\hfill

{\bf Proof:}
Square of Hermitian norm of $\eta$ is equal to $\int \eta\wedge\eta\wedge\omega^{1,1}$, but in dimension 3
we have the equality $\eta\wedge\eta\wedge\omega^{1,1}=\eta\wedge\eta\wedge\omega$; the latter form
 is exact,
therefore $\eta=0$.
\endproof

\hfill

\corollary\label{ddc11}
Let $M$ be a compact Hermitian symplectic threefold, $\alpha$ is a $d$-closed, $d^c$-exact $(1,1)$-form.
Then $\alpha=dd^cf$ for some function $f$.
\endproof

\hfill

%%%%%%%%%%%%%%%%%%%%%%%%%%%%%%%%%%%%%%%%%%%%%%%%%%%%%%%%%%%%%%%%%%%%%%%
\section{Applications}
%%%%%%%%%%%%%%%%%%%%%%%%%%%%%%%%%%%%%%%%%%%%%%%%%%%%%%%%%%%%%%%%%%%%%%%

%It is rather well-known and easy to show that every bounded bicomplex over a field
%decomposes into a direct sum of subbicomplexes of the following kind:
%
%(TO DEPICT)
%
%Note that $dd^c$-lemma is equivalent to the absence of all zig-zags, and
%degeneracy of Fr\"olicher spectral sequence is equivalent to the 
%absence of zig-zags of odd length.
%
%By \ref{holoforms}, \ref{ddc11} and degeneracy of Fr\"olicher spectral
%sequence proved in \cite{Ca} possible zig-zags in the Dolbeault 
%bicomplex are of the following kind:
%
%\hfill

\theorem (Gauduchon, \cite{Gauduchon}).
For a complex manifold $M$, $dd^c$-lemma for $(1,1)$-forms is equivalent
to the equality $b^1=2h^{1,0}$.

\hfill

{\bf Proof:} Consider the cohomology 
sequence associated to the short exact sequence of sheaves of the form
$0 \arrow \sqrt{-1}\R \arrow \mathcal{O} \stackrel{Re}{\arrow} \H \arrow 0$, where
$\H$ is the sheaf of pluriharmonic functions.
The relevant piece looks like $$... \stackrel{0}{\arrow} H^1(M,\sqrt{-1}\R) \arrow H^1(M,\mathcal{O})
\arrow H^1(M,\H) \arrow H^2(M,\sqrt{-1}\R) \arrow ...$$
It is well-known that 

$$H^1(M,\H)=\frac{\Ker d:\A^{1,1} \arrow \A^3}{\Im dd^c: \A^0 \arrow \A^{1,1}}.$$

So $dd^c$-lemma for $(1,1)$-forms holds if and only if the third arrow is an
isomorphism, and, by exactness, if and only if the first arrow is an isomorphism. \endproof

\hfill

So, the equality $b^1=2h^{1,0}$ holds on compact Hermitian symplectic threefolds.
 
\hfill

It follows that we have the Hodge decomposition on the first cohomology of $M$:
$H^1(M,\C)=H^{0,1}(M) \oplus H^{1,0}(M)$, and $\dim H^{0,1}=\dim H^{1,0}$.
So the rank of the abelian group $H^1(M,\Z)$ is equal to the
dimension of the real vector space $H^{0,1}(M)$. It follows that
the Albanese torus is defined correctly and we have the Albanese map
$\Alb: M \arrow H^{0,1}(M)^*/H_1(M,\Z)$. Its image $\Alb(M)$ is a subvariety 
(possibly singular) of 
a torus.

\hfill

\proposition
Suppose $\dim \Alb(M)=3$. Then $M$ is K\"ahler.

\hfill

{\bf Proof:}
If $\Alb(M)$ is smooth, then $\Alb$ is an immersion, and pullback of
the K\"ahler form $\Alb^*\omega$ is the K\"ahler form on $M$.
Otherwise, we can desingularize the morphism $\Alb$ to obtain the K\"ahler
metric on some manifold $\tilde{M}$ bimeromorphic to $M$ ($M$ is then a manifold
in the Fujiki class C). On the other hand, $M$ admits an SKT structure (\ref{SKT}).
From the theorem of Chiose (\cite{Chiose}) it follows that $M$ is K\"ahler. \endproof

\hfill

\remark
It would be interesting to know what one can extract from the Albanese map if 
$\dim \Alb(M)=1$ or $2$. For example, if $\Alb(M)$ is a smooth curve $C$, fibers
of $\Alb(M)$ are Hermitian symplectic (and therefore K\"ahler) surfaces,
and the pullback of the volume form $\Alb^* \Vol_C$ is a closed, non-exact $(1,1)$-form
on $M$. By $dd^c$-lemma for $(1,1)$-forms and \ref{holoforms}, it could not be cohomologous
to a form of type $(2,0)+(0,2)$. By a theorem of Michelsohn (\cite{Michelsohn}), in that
situation there exists a {\em balanced} metric on $M$, that is, a Hermitian form $\omega$
such that $d\omega^{\dim M - 1}=0$. Actually, the smoothness of $C$ is not necessary,
because a manifold bimeromorphic to a balanced manifold is balanced itself (\cite{AB}).

\addcontentsline{toc}{section}{References}

\noindent {\sc {\bf G.P.}:
{\sc Laboratory of Algebraic Geometry,\\
National Research University HSE,\\
Department of Mathematics, 7 Vavilova Str. Moscow, Russia,}\\
\tt  datel@mail.ru}.
\end{document}